\documentclass[final,leqno,letterpaper]{article}
\usepackage{amsmath}
\usepackage{hyperref}
\usepackage{graphicx}
\usepackage{authblk}
\usepackage{algpseudocode}
\usepackage{amsthm}

\hypersetup{
    pdftitle={An Improved Formula for Jacobi Rotations},
    pdfauthor={Carlos F. Borges},
    pdfkeywords={Rotation, Stable}
    }

\newtheorem{algorithm}{Algorithm}

\def\bne{\begin{equation}}
\def\ene{\end{equation}}
\def\bmat{\left[ \begin{array}}
\def\endbmat{\end{array} \right]}

\title{An Improved Formula for Jacobi Rotations} 
\date{June 30, 2017}

\author{Carlos F. Borges}
\affil{Department of Applied Mathematics\\Naval Postgraduate School\\Monterey CA 93943\\Email: borges@nps.edu}

\begin{document}
\maketitle

\begin{abstract}
We present an improved form of the algorithm for constructing Jacobi rotations. This is simultaneously a more accurate code for finding the eigenvalues and eigenvectors of a real symmetric $2 \times 2$ matrix.
\end{abstract}

\pagestyle{myheadings}
\thispagestyle{plain}

Given a $2 \times 2$ real symmetric matrix 
$$
A = \bmat{c c}
a_{pp} & a_{pq} \\
a_{qp} & a_{qq}
\endbmat
$$
the standard {\em stable} algorithm for constructing a Jacobi rotation that diagonalizes $A$ so that
$$
\bmat{r r} c & s \\ -s & c \endbmat^T
\bmat{c c} a_{pp} & a_{pq} \\ a_{qp} & a_{qq} \endbmat
\bmat{r r} c & s \\ -s & c \endbmat
= \bmat{c c} \lambda_1 & 0 \\ 0 & \lambda_2 \endbmat
$$
where the $\lambda_i$ are unordered eigenvalues of $A$, can be found in \cite{gvl} as well as a number of online resources. If we include the stable computation of the eigenvalues the algorithm can be coded is as follows:

\vspace{.2in}
\begin{algorithm} The Standard Approach

\hrule
\begin{algorithmic}
\If{$a_{pq} \neq 0$}
	\State $\delta \gets (a_{qq} - a_{pp})/(2a_{pq})$
	\If{$\delta \geq 0$}
		\State $t \gets 1/(\delta+\sqrt{1+\delta^2})$
	\Else
		 \State $t \gets 1/(\delta-\sqrt{1+\delta^2})$
	\EndIf
\Else 
	\State $t \gets 0$
\EndIf

\State $c \gets 1/\sqrt{1+t^2}$
\State $s \gets tc$
\State $\lambda_1 = a_{pp}-ta_{pq}$
\State $\lambda_2 = a_{qq}+ta_{pq}$

\end{algorithmic}
\hrule
\end{algorithm}
\vspace{.2in}

The algorithm above was constructed to avoid unecessary overflow that might occur in an interim calculation but this approach has been superseded as most modern computing environments provide a function called {\tt hypot(a,b)} that can better deal with this issue\footnote{The math library function {\tt hypot(a,b)} calculates $\sqrt{a^2+b^2}$ in a manner that avoids unecessary overflow or underflow when the arguments are badly scaled.} and we can improve numerical performance if we take advantage of it. A better form of this algorithm is as follows:

\vspace{.2in}
\begin{algorithm} The Improved Approach

\hrule
\begin{algorithmic}
\If{$a_{pq} \neq 0$}
	\State $\delta \gets (a_{qq} - a_{pp})/2$
	\If{$\delta \geq 0$}
		\State $t \gets a_{pq}/(\delta+{\tt hypot}(a_{pq},\delta))$
	\Else
		 \State $t \gets a_{pq}/(\delta-{\tt hypot}(a_{pq},\delta))$
	\EndIf
\Else 
	\State $t \gets 0$
\EndIf

\State $c \gets 1/\sqrt{1+t^2}$
\State $s \gets tc$
\State $\lambda_1 = a_{pp}-ta_{pq}$
\State $\lambda_2 = a_{qq}+ta_{pq}$
\end{algorithmic}
\hrule
\end{algorithm}
\vspace{.2in}

In order to compare the performance of the two approaches we will see how well the Jacobi rotation (which is in essence the matrix of normalized eigenvectors $V$) and the computed eigenvalues satisfies the fundamental identity $AV = V\Lambda$. 

For our testing we will look at three algorithms, the two described in this note as well as the appropriate algorithm from LAPACK for solving the real symmetric eigenvalue problem. The LAPACK codes \cite{LUG} are the 'industry standard' but we should note that they are designed for arbitrary sized real symmetric matrices, and not specifically for a $2\times2$. Our tests will proceed by comparing the magnitude of $\Vert AV-V\Lambda \Vert_F$ for each of the three algorithms over a set of test matrices. We describe the test as it is implemented to see how the algorithms perform over a range of 'scales' for the off-diagonal element $a_{pq}$: 
\begin{enumerate}
\item
Generate a random set of 100,000 real symmetric matrices with elements distributed according to a standard normal distribution.
\item
Generate a range of variances over which $a_{pq}$ will be scaled for the test.
\item
For each variance value multiply the off-diagonal element of every matrix in the test set by the square root of that value so that the $a_{pq}$ values for the test set exhibit the proper variance.
\item
Use each of the three algorithms to find $V$ and $\Lambda$ for every matrix in the scaled test set and then compute the average value of $\Vert AV-V\Lambda \Vert_F$ for each.
\item
Plot the output.
\end{enumerate}

A similar approach is used for running $a_{pp}$ through a range of values (there is no need to do so for $a_{qq}$ as that would yield the same result). All testing was done in Matlab 2016a on an Intel(R) Core(TM) i7-2600 CPU.

The results of all of these tests appear in graphs at the end of this article and they demonstrate the superiority of the improved algorithm since it functions at least as well as the best of the other two in general and is better than both of the others in cases of extreme scaling.

\begin{figure}[ ht ]
\begin{center}
\includegraphics[width=6in]{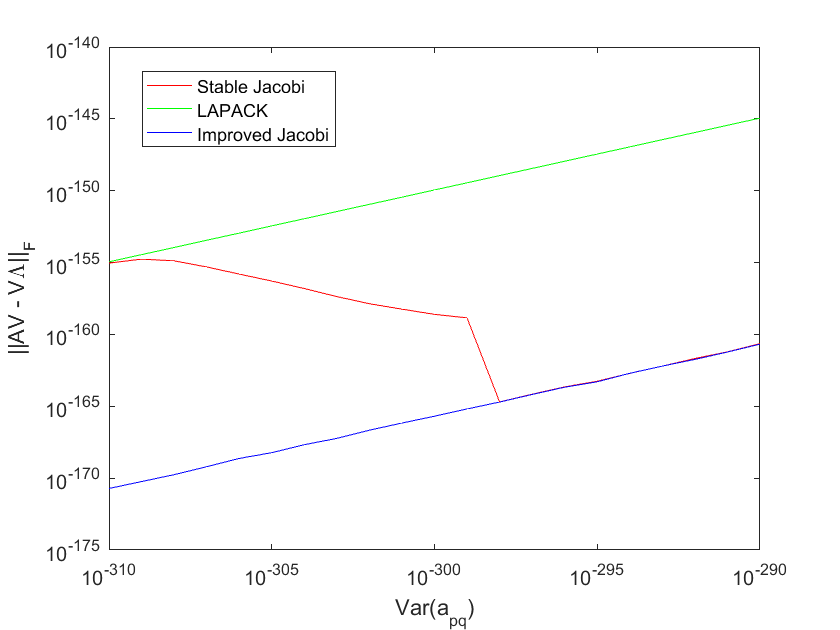}
\caption{Comparison of the three algorithms when the variance (i.e. scale) of $a_{pq}$ is manipulated. Test matrix elements distributed ${\cal N}(0,1)$ except for $a_{pq}$ which is a zero mean normal with the variance ranging over the values shown on the x-axis.}
\end{center}
\end{figure}

\begin{figure}[ ht ]
\begin{center}
\includegraphics[width=6in]{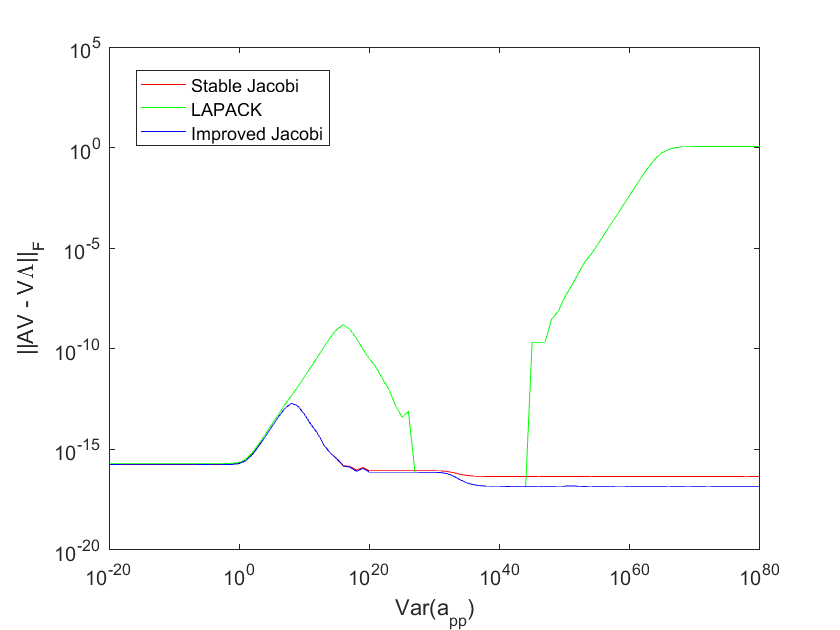}
\caption{Comparison of the three algorithms when the variance (i.e. scale) of $a_{pp}$ is manipulated to be large. Test matrix elements distributed ${\cal N}(0,1)$ except for $a_{pp}$ which is a zero mean normal with the variance ranging over the values shown on the x-axis.}
\end{center}
\end{figure}

\begin{figure}[ ht ]
\begin{center}
\includegraphics[width=6in]{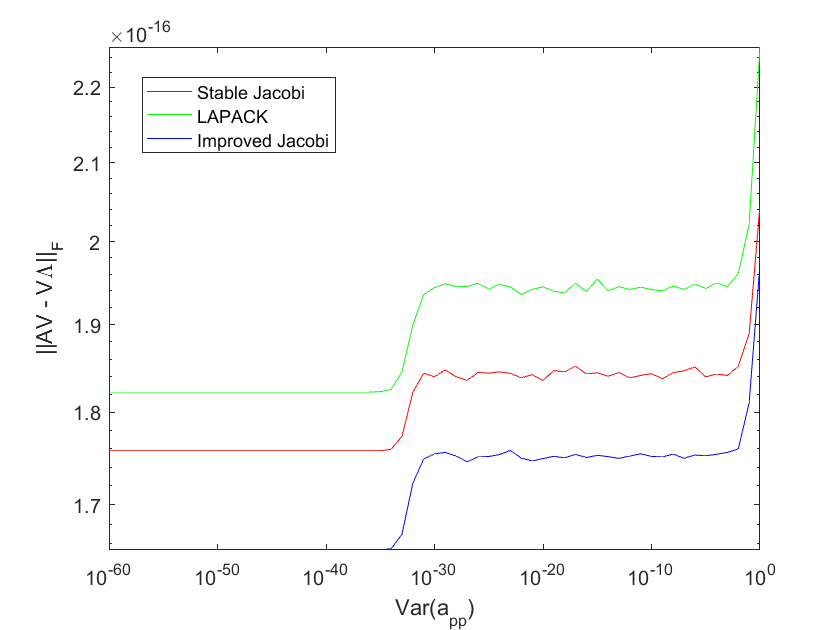}
\caption{Comparison of the three algorithms when the variance (i.e. scale) of $a_{pp}$ is manipulated to be small. Test matrix elements distributed ${\cal N}(0,1)$ except for $a_{pp}$ which is a zero mean normal with the variance ranging over the values shown on the x-axis.}
\end{center}
\end{figure}

\end{document}